\documentclass[11pt,a4paper]{amsart}
\usepackage{amsfonts}
\usepackage{amsthm}
\usepackage{amsmath}
\usepackage{amssymb}
\usepackage{amscd}
\usepackage{t1enc}
\usepackage[margin=2.9cm]{geometry}
\usepackage{epstopdf} 
\usepackage{times}
\usepackage{enumerate}
\usepackage{changepage}
\usepackage{mfirstuc}  
\usepackage{tikz}
\usepackage{cite}
\usepackage{booktabs}
\usepackage{tabularx}
\numberwithin{equation}{section}

\theoremstyle{plain}
\newtheorem{theorem}{Theorem}[section]
\newtheorem{lemma}[theorem]{Lemma}

\newtheorem{cor}[theorem]{Corollary}
\newtheorem{proposition}[theorem]{Proposition}
\newtheorem{problem}[theorem]{Problem}
\newtheorem{question}[theorem]{Question}
\newtheorem{definition}[theorem]{Definition}
\newtheorem{conjecture}[theorem]{Conjecture}
\theoremstyle{remark}
\newtheorem{remark}[theorem]{Remark}
\newtheorem{example}[theorem]{Example}

\newcommand{\nfp}{\mathbb{F}_{p}}
\newcommand{\nfq}{\mathbb{F}_{q}}

\newcommand{\F}{\mathbb{F}}
\newcommand{\C}{\mathbb{C}}

\newcommand{\Nm}{\operatorname{N}}
\newcommand{\Tr}{\operatorname{Tr}}

\usepackage{xcolor}
\definecolor{eyecolor}{RGB}{204,232,207}
\pagecolor{eyecolor!90}

\begin{document}

\title{Higher Kloosterman sums in finite semi-simple algebras}

\author{Daqing Wan}
\address{Center for Discrete Mathematics, Chongqing University, Chongqing 401331, China}
\address{College of Mathematics and Statistics, Chongqing University, Chongqing 401331, China}

 \email{dwan@math.uci.edu} 
  


\begin{abstract} An asymptotic formula with a square root error term is obtained for higher Kloosterman sums over semi-simple 
	algebras over a finite field. This asymptotic formula has an intermediate 
	main term whose vanishing depends on the vanishing of a one variable exponential sum over the finite field. 
	This leads to a non-vanishing conjecture for such exponential sums, which can be viewed as an exponential 
	sum analogue of Lenstra's theorem on the Carlitz-Wan conjecture for permutation polynomials.

	\end{abstract}

\maketitle
\setcounter{tocdepth}{1}

\section{Introduction}

Let $\F_q$ be the finite field of $q$ elements with characteristic $p$. By the Artin-Wedderburn theorem, 
any finite  semi-simple algebra $M$ over $\F_q$ of degree $n\geq 2$ is isomorphic to a finite product of matrix algebras over $\F_q$. 
Namely, 
\begin{equation}\label{B}
M = M_{d_1}(\F_{q^{n_1}}) \times \cdots \times M_{d_k}(\F_{q^{n_k}}), \ \ n_1d_1^2+\cdots + n_kd_k^2 = n,\ |M| =q^n, 
\end{equation}
where $\F_{q^{n_i}}$ denotes the finite field extension of $\F_q$ with degree $n_i$, and $M_{d_i}(\F_{q^{n_i}})$ denotes the $d_i\times d_i$ 
matrix algebra over $\F_{q^{n_i}}$. 
This factorization of $M$ is uniquely determined  up to permutations. The group of invertible elements in $M$ is 
\begin{equation}
M^* = {\rm GL}_{d_1}(\F_{q^{n_1}}) \times \cdots \times {\rm GL}_{d_k}(\F_{q^{n_k}}), \ \ |M^*| = \prod_{i=1}^k |{\rm GL}_{d_i}(\F_{q^{n_i}})|.
\end{equation}
The semi-simple algebra $M$ is called  \'{e}tale (or commutative) if all $d_i=1$, in which case, $M$ is denoted by $B$.  
For $x \in M$, multiplication by $x=(A_1, \cdots, A_k)$ on the left acts on the reduced $\F_q$-vector space 
$$R_M: = \bigoplus_{i=1}^k \F_{q^{n_i}}^{d_i}, \ \ \ \ \dim_{\F_q}(R_M) = \sum_{i=1}^k d_in_i.$$
Its trace ${\rm Tr}_M(x)$ and norm ${\rm N}_M(x)$, called reduced trace 
and reduced norm in the literature, are well defined elements in $\F_q$. 

\begin{definition}
Let $\psi: \F_q \longrightarrow \C^*$ be a fixed non-trivial additive character. Let $f(t) \in \F_q[t]$ be a polynomial of degree $r\geq 1$. 
For $b\in \F_q^*$,  the higher Kloosterman sum over $M$ associated to $f(t)$ is defined by
$$K_{M,f}(b, \psi): = \sum_{x \in M^*, \ {\rm N}_M(x)=b} \psi(f({\rm Tr}_M(x))).$$
In the special case $f(t)=t$, this sum $K_{M,t}(b, \psi)$ is called the usual Kloosterman sum over $M$,  denoted by 
$$K_{M}(b, \psi): = \sum_{x \in M^*, \ {\rm N}_M(x)=b} \psi({\rm Tr}_M(x)).$$
\end{definition}

In a previous note \cite[Theorem 1.6]{Wa2}, the following square root estimate is obtained for the usual Kloosterman sum $K_M(b, \psi)$ over $M$: 
\begin{align}\label{klsum}
\left| K_{M}(b, \psi): = \sum_{x \in M^*, \ {\rm N}_M(x)=b} \psi({\rm Tr}_M(x))\right| \leq  (\sum_{i=1}^k d_in_i )q^{\frac{n-1}{2}} \leq n q^{\frac{n-1}{2}}.
\end{align}
This estimate extends Deligne's classical estimate \cite[Equation 7.2.5]{De} from commutative $M=B$ (when all $d_i=1$) to general non-commutative $M$.  
Note that the constant $\sum_{i=1}^k d_in_i$ is strictly smaller than $n$ when $M$ is non-commutative, i.e., when some $d_i>1$.

Motivated by a number of applications, we now study the higher Kloosterman sum $K_{M,f}(b, \psi)$ when $r= \deg(f)>1$. 
In this higher degree setting, new phenomenon occurs as there will be a non-trivial intermediate main term. Thus, it is not possible to obtain a simple 
square root estimate as the case $r=1$. 
The best one can hope for is an asymptotic formula with both an intermediate main term and a square root error term. 
If $r$ is divisible by $p$, a standard Artin-Schreier reduction reduces $f(t)$ to a lower degree polynomial. Thus, we can always assume that $(r, p)=1$ as the constant polynomial case is trivial. 
Our main result is 

\begin{theorem}\label{thm}
Let $M$ be a finite semi-simple algebra over $\F_q$ of degree $n=\sum_{i=1}^k d_i^2 n_i \geq 2$ as given in equation (\ref{B}) with $m=\sum_{i=1}^k d_in_i$.   
Let $f(t) \in \F_q[t]$ be a polynomial of degree $r\geq 1$ with $(r, p)=1$. 
For $b\in\F_q^*$, we have 
	
\begin{align}\label{kl5}
\left|K_{M, f}(b, \psi) - C_M\sum_{a\in \F_q}\psi(f(a)) \right|   \leq\left( rm -r +1\right) q^{\frac{n-1}{2}}. 
		\end{align}
		where $C_M$ is the explicit coefficient 
		\[ C_M: = \frac{|M^*| +(-1)^{\sum_i d_i}(q-1)q^{\frac{n-m}{2}}}{q(q-1)}. \]
		\end{theorem}	
In the case $f(t)=t$, the main term becomes zero and estimate (\ref{kl5}) reduces to estimate (\ref{klsum}). 		
Note that the constant $rm-r+1$ is at most $rn$, and is strictly smaller than $rn$ if either $r>1$ or $M$ is non-commutative. 
In proving Theorem \ref{thm}, our main idea is to apply the reduction in \cite{Wa2} to reduce the sum $K_{M, f}(b, \psi)$ to 
an algebraic exponential sum, and then apply the theory of $\ell$-adic cohomology. In addition to Deligne's Kloosterman sheaf over $\mathbb{G}_m$, we also need Katz's theory \cite{K} of Airy sheaf over $\mathbb A^1$ which plays a crucial role.  In the end, the problem is reduced to calculating the 
cohomology of the tensor product of these two sheaves. 

The existence of a non-trivial main term in (\ref{kl5}) means that there are non-trivial cohomology beyond the middle dimension. 
The main term is clearly zero if $r=\deg(f)=1$, or more generally if $f(t)$ is a permutation polynomial over $\F_q$, that is, if $f$ induces a 
permutation of $\F_q$. 		
		
\begin{cor}\label{cor4} 
Let $M$ be a finite semi-simple algebra over $\F_q$ of degree $n=\sum_{i=1}^k d_i^2 n_i \geq 2$ as given in equation (\ref{B}) with $m=\sum_{i=1}^k d_in_i$.   
 Let $f(t) \in \F_q[t]$ be a permutation polynomial of degree $r\geq 1$ with $(r, p)=1$.  For $b\in\F_q^*$, we have 
	\begin{align}
		\left| K_{M, f}(b, \psi) \right| 
		\leq \left( rm -r +1 \right)q^{\frac{n-1}{2}}.
	\end{align}

\end{cor}

\begin{remark}
Without the assumption that $f(t)$ is a permutation polynomial, Theorem (\ref{thm}) gives the following much weaker Lang-Weil type bound 
\begin{align}
		\left| K_{M, f}(b, \psi) \right| 
		\leq (r-1)|C_M| \sqrt{q} + (rm-r+1) q^{\frac{n-1}{2}}.
	\end{align}
	This bound is essentially best possible in general, as in some super-singular cases, one could have 
	$$|\sum_{a\in \F_q}\psi(f(a))|=(r-1)\sqrt{q},$$
	attaining the Weil bound. 
	\end{remark}

\begin{remark}	
In the special case when 
$M=M_{d}(\F_q)$ and $f(t)=t^r$, the sum $K_{M, t^r}(b, \psi)$ is already 
studied in Kim \cite[Theorem 5.1]{Ki} who applied 
the Bruhat decomposition to recursively reduce the sum to the study of the following higher Kloosterman sum over the geometric $M= \F_q^n$:  
\begin{align}
		K_{\F_q^n, t^r}(b, \psi): =  \sum_{x_i \in \F_q^*, \ x_1\cdots x_n=b} \psi\big((x_1+\cdots +x_n)^r\big), \ \ b \in \F_q^*.
	\end{align}
But no estimate for $K_{\F_q^n, t^r}(b, \psi)$ is given in \cite{Ki}. 
Corollay \ref{cor4} shows that for $(r, q-1)=1$, one obtains the following square root estimate:
	\begin{align}
		\left| K_{\F_q^n, t^r}(b, \psi): =  \sum_{x_i \in \F_q^*, \ x_1\cdots x_n=b} \psi\big((x_1+\cdots +x_n)^r\big) \right| 
		\leq (nr-r+1) q^{\frac{n-1}{2}} , \ \ b \in \F_q^*.
	\end{align}	
 For the connection of Theorem \ref{thm} to hypergeometric sums over reductive groups 
studied abstractly in Fu-Li \cite{FL}, see Section $4$ for additional discussions. A number of open problems including an exotic version of 
the higher Kloosterman sum will be discussed there as well.   
\end{remark}

To conclude this introduction, we give a brief discussion about the vanishing and non-vanishing of 
the main term in Theorem \ref{thm}.  		
First, by the Weil bound, the sum $\sum_{a\in \F_q}\psi(f(a))$ is bounded by $(r-1)\sqrt{q}$.  
In concrete applications, to obtain the best square root cancellation 
for $K_{M, f}(b, \psi)$, we want to know when the main term is zero. This raises the following question. 

\begin{question}\label{Qu}
Let $f(t) \in \F_q[t]$ be a polynomial of degree $r\geq 1$ with $(r, p)=1$. Decide when $\sum_{a\in \F_q}\psi(f(a))=0$ and when 
$\sum_{a\in \F_{q^e}}\psi ({\rm Tr}_e(f(a)))=0$ for infinitely many $e$? Here ${\rm Tr}_e$ denotes the trace map from $\F_{q^e}$ 
to $\F_q$. 
\end{question}

Recall that a polynomial $f(t)\in \F_q[t]$ is called exceptional over $\F_q$ if $f(t)$ is a permutation polynomial over $\F_{q^e}$ for 
infinitely many $e$. For large $q$, say $q> r^4$, the Weil bound and a theorem of Cohen \cite{Co} show that $f(t)\in \F_q[t]$ is an exceptional polynomial over $\F_q$ if and only 
if $f(t)$ is a permutation polynomial over $\F_q$. 
Question \ref{Qu} is thus related to permutation polynomials and exceptional polynomials over $\F_q$, but could be somewhat harder.  
For instance, no polynomial time algorithm seems to be known for Question \ref{Qu}.  
In contrast, there is a polynomial time algorithm for testing if $f(x)$ is a permutation polynomial over $\F_q$, see Kayal \cite{Kay}. 
The sum $\sum_{a\in \F_q}\psi(f(a))$ is zero if 
$f(t)$ is a permutation polynomial over $\F_q$. The converse is only true if $q=p$, but not true if $q$ is a proper power of $p$. 
If $f(t)=t^r$ with $(r, q-1)=1$ or if $f(t)$ is 
a Dickson polynomial of degree $r$ with $(r, q^2-1)=1$, then $f(t)$ is a permutation polynomial and also an exceptional polynomial over $\F_q$. 
For $f(t)=t^r$, a congruence argument modulo $(r, q-1)$ shows that 
$\sum_{a\in \F_q}\psi(a^r)=0$ if and only if $(r, q-1)=1$. For Dickson polynomials of degree $r$, we do not know 
a simple necessary and sufficient condition for $\sum_{a\in \F_q}\psi(f(a))=0$. The vanishing exponential sum problem is already quite subtle for binomial $f(x)=x^r +bx$. 
In this case, we mention an old conjecture of Helleseth \cite{He} (1976) which says that for 
$(r, q-1)=1$, $r \equiv 1 \mod (p-1)$, there is at least one $b \in \F_q^*$ such that $\sum_{a\in \F_q}\psi(a^r +ba)=0$ for $q>2$. 
This conjecture seems still wide open.

Another related problem is the Carlitz-Wan conjecture \cite{Wa1} which says that if $(r, q-1)>1$, then there are no exceptional polynomials of degree $r$ over $\F_q$. 
In the special case that $q$ is odd and $r$ is even, this is the original Carlitz conjecture which was first proved by Fried-Guralnick-Saxl \cite{FGS} using a rather sophisticated 
method involving the classification of finite simple groups. The general Carlitz-Wan conjecture was first proved by Lenstra using a much simpler algebraic 
number theory method, see Cohen-Fried \cite{CF}, Guralnick-M\"uller \cite{GM}, Ding-Xiong-Zhang \cite{DXZ}  and Hu-Zhang \cite{HZ} 
for various proofs, generalizations and historical remarks.  
The following non-vanishing conjecture can be viewed as an exponential sum analogue of Lenstra's theorem.  

\begin{conjecture}\label{con} Let $f(t) \in \F_q[t]$ be a polynomial of degree $r\geq 1$ with $(r, p)=1$. Assume $(r, q-1)>1$. 
If the characteristic $p$ is sufficiently large compared to $r$,  then $\sum_{a\in \F_q}\psi(f(a)) \not =0$. 
\end{conjecture}

This conjecture has a positive answer if $q=p$,  in which case $\sum_{a\in \F_p}\psi(f(a))=0$ if and only if $f(t)$ is a permutation 
polynomial over $\F_p$. One can then apply the results on permutation polynomials. If $q$ is a proper power of $p$, the above 
conjecture seems to be significantly deeper, see the last section for some evidence and a heuristic argument supporting the 
conjecture. As we shall see,  the condition that the characteristic $p$ is sufficiently large compared to $r$ cannot be replaced by the 
weaker condition that $q$ is sufficiently large compared to $r$. 

The rationality of the L-function shows that the following sequence of exponential sums 
$$S_{q^e}(f): =\sum_{a\in \F_{q^e}}\psi ({\rm Tr}_e(f(a))), \ e=1,2 , \cdots$$
is a sequence of algebraic integers satisfying a linear recurrence relation. 
The Skolem-Mahler-Lech theorem then implies that its zero set 
$$Z_q(f): =\{ e \in \mathbb{N} \ | \   S_{q^e}(f)=0 \}$$ 
is the union of a finite set and finitely many arithmetic progressions. 
The celebrated Skolem problem in this case is to decide if the zero set $Z_q(f)$ is empty.  
For a general algebraic linear recurrence sequence, the Skolem problem is at least ${\rm NP}$-hard, and  currently not even known to be decidable, 
except for algebraic linear recurrence sequences of order at most $4$, see Bacik \cite{Ba}. Note that the weaker question of deciding if the zero set $Z_q(f)$ is infinite is decidable. In our special exponential sum case, we believe that the Skolem problem should be decidable and perhaps even has a polynomial 
time algorithm. From this point of view, Conjecture \ref{con} 
says that if $(r, q-1)>1$ and the characteristic $p$ is sufficiently large compared to the degree $r$ of $f(t)$,  then 
the zero set $Z_q(f)$ is empty.

\begin{remark} This paper was written in 2025 before the AI age, and it is presented here as is, except for this remark and the updated references. 
Now, in the AI age, it is tempting to use the AI tool to see what can be proved toward Conjecture \ref{con}. A few rounds of discussions with GPT-6 Astra 
lead to the following partial results proved using $\ell$-adic methods and $p$-adic methods.  
\begin{itemize}
\item Conjecture \ref{con} holds if  $q=p^e$ and $p> r^{4e}$. 
\item Conjecture \ref{con} holds if $q=p^e$ and $f(t) \in \mathbb{F}_p[t]$. 
\item  Conjecture \ref{con} holds if the Galois group of $f(t)-x$ over $\bar{\mathbb{F}}_q(x)$ is doubly transitive. 
\end{itemize}
These AI assisted partial results suggest that Conjecture \ref{con} might be sufficiently non-trivial. 
\end{remark}

\section{Points on the norm-trace equations}

For convenience, we first recall the reduction \cite{Wa2} for counting the number of elements in $M^*$ with fixed trace and norm. 
This reduction reduces our problem to estimating certain algebraic exponential sums. 
It will be our starting point to prove Theorem \ref{thm}. 

Let $ \Tr_{n_i}=\Tr_{\F_q^{n_i} / \F_q}$ be the trace map from $\F_{q^{n_i}}$ to $\F_q$. Similarly, let $\Nm_{n_i}=\operatorname{Norm}_{\F_{q^{n_i}} / \F_q}$ be the norm map from $\F_{q^{n_i}}$ to $\F_q$. For $a \in \F_q$ and $b \in \F_q^*$, let 
\begin{align*}
	\Nm_{M}(a,b)=
	\#\left\{ (x_1,\cdots, x_k)\ \left| \ \ x_i\in {\rm GL}_{d_i}(\F_{q^{n_i}}), \ \ 
	\sum_{i=1}^{k}\Tr_{n_i}({\rm tr}(x_i))=a,\ \ \prod_{i=1}^{k}\Nm_{n_i}({\rm det}(x_i))=b
	\right. \right\},
\end{align*}
where ${\rm tr}(x_i)$ and ${\rm det}(x_i)$ denote the trace and determinant of the $d_i\times d_i$ matrix $x_i$. It is clear that 
\begin{align*}
	\Nm_{M}(a,b)=
	\#\left\{ x \in M \ \left| \  \Tr_{M}(x)=a,\ \ \Nm_{M}(x)=b
	\right. \right\} = \#\left\{ x \in M^* \ \left| \  \Tr_{M}(x)=a,\ \ \Nm_{M}(x)=b
	\right. \right\} .
\end{align*}

Let $\psi$ be a fixed non-trivial additive character of $\F_q$. 
For a multiplicative character of $\F_q^*$,  let $G(\chi)$ denote the Gauss sum
$$
G(\chi)=\sum_{x\in\F_q^*}\psi(x)\chi(x).
$$
Note that $G(\chi) =-1$ if $\chi$ is the trivial character, and $|G(\chi)|=\sqrt{q}$ if $\chi$ is non-trivial.  
Let $\psi_{n_i}=\psi\circ {\Tr_{n_i}}$ be the lifted additive character of
$\F_{q^{n_i}}$, and let $\chi_{n_i} = \chi \circ {{\rm N}_{n_i}}$ be the lifted multiplicative character of $\F_{q^{n_i}}^*$. The Hasse-Davenport 
relation gives 
$$G(\chi_{n_i}, \psi_{n_i} ) =  \sum_{x\in\F_{q^{n_i}}^*}\psi_{n_i}(x)\chi_{n_i}(x) = (-1)^{n_i-1}G(\chi)^{n_i}.$$
Now, we relate $\Nm_{M}(a,b)$ to algebraic exponential sums over $\F_q$. 
 Let $m=\sum_i d_in_i$. 
Let $x_i$ below run over ${\rm GL}_{d_i}(\F_{q^{n_i}})$, and let $\chi$ run over all multiplicative characters of $\F_q^*$. By orthogonality of characters, one computes 
\begin{align}\label{eq1}
	q(q-1) \Nm_M(a, b)
	&=\sum_{\left(x_1, \cdots, x_k\right)} \sum_{v \in \F_q} \psi\left(v(\sum_{i=1}^{k}\Tr_{n_i}\left({\rm tr}(x_i)\right)-a)\right) 
	\sum_{\chi} \chi\left(b^{-1} \prod_{i=1}^{k}\Nm_{n_i}\left({\rm det}(x_i)\right)\right) \nonumber\\
	& =|M^*| + \sum_{v\neq 0} \psi(-a v) \sum_{\chi} \bar{\chi}(b) \sum_{\left(x_1, \cdots, x_k\right)} \prod_{i=1}^{k}\psi_{n_i}\left({\rm tr}(vx_i)\right)
	\prod_{i=1}^{k}\chi_{n_i}\left({\rm det}(x_i)\right) \nonumber\\
	& =|M^*| +\sum_{v \neq 0} \psi(-a v) \sum_{\chi} \bar{\chi}(bv^{\sum_i d_in_i})\prod_{i=1}^{k}\left(\sum_{x_i} \psi_{n_i}\left({\rm tr}(x_i)\right) \chi_{n_i}\left({\rm det}(x_i)\right)\right) \nonumber\\
&=|M^*| +\sum_{v \neq 0} \psi(-a v) \sum_{\chi} \bar{\chi}(bv^{m})\prod_{i=1}^{k} G_{{\rm GL}_{d_i}(\F_{q^{n_i}})}(\chi_{n_i}, \psi_{n_i}), 
\end{align}
where the last factor is the Gauss sum over the general linear group ${\rm GL}_{d_i}(\F_{q^{n_i}})$ associated to the 
character $(\chi_{n_i}, \psi_{n_i})$ of $\F_{q^{n_i}}^* \times \F_{q^{n_i}}$, defined by
$$G_{{\rm GL}_{d_i}(\F_{q^{n_i}})}(\chi_{n_i}, \psi_{n_i}) = \sum_{x_i\in {\rm GL}_{d_i}(\F_{q^{n_i})} }\psi_{n_i}\left({\rm tr}(x_i)\right) \chi_{n_i}\left({\rm det}(x_i)\right).$$
By Eichler's formula \cite{EI} (see also \cite{La}\cite{Ki} for new proofs),  
 we have
$$G_{{\rm GL}_{d_i}(\F_{q^{n_i}})}(\chi_{n_i}, \psi_{n_i}) =q^{n_i {d_i \choose 2}} G(\chi_{n_i}, \psi_{n_i})^{d_i}.$$
This together with the Hasse-Davenport relation gives 
$$G_{{\rm GL}_{d_i}(\F_{q^{n_i}})}(\chi_{n_i}, \psi_{n_i}) 
 =(-1)^{d_i(n_i-1)}q^{n_i {d_i \choose 2}} G(\chi)^{d_in_i}.$$
Recall that $m=\sum_i d_in_i$ and $n= \sum_i d_i^2n_i$. 
It follows that we have 
\begin{align}
q(q-1) \Nm_M(a, b)&=|M^*| +(-1)^{\sum_i d_i(n_i-1)}q^{\sum_i n_i {d_i \choose 2}}\sum_{v \neq 0} \psi(-a v) \sum_{\chi} \bar{\chi}(bv^{m})
	G(\chi)^{\sum_i d_in_i} \nonumber\\
 &=|M^*| +(-1)^{m-\sum_i d_i} q^{\frac{n-m}{2}} S_m(a, b),  	
		\end{align}
		where 
\begin{align}\label{sum}
S_m(a, b)&= \sum_{v \neq 0} \psi(-a v) \sum_{\chi} \bar{\chi}(bv^{m})
	G(\chi)^{m}.
	\end{align}
This gives the following reduction 
\begin{lemma}\label{lem} We have the relation
\begin{align}
\Nm_M(a, b)&=\frac{|M^*|}{q(q-1)} + \frac{(-1)^{m-\sum_i d_i} q^{\frac{n-m}{2}}}{q(q-1)} S_m(a, b).  	
		\end{align}
\end{lemma}
This key lemma was used to prove the following estimate \cite[Theorem 1.1]{Wa2} for $\Nm_M(a, b)$.  

\begin{theorem}\label{thm2} Let $M$ be a finite semi-simple algebra over $\F_q$ of degree $n=\sum_{i=1}^k d_i^2 n_i \geq 2$ as given in equation (\ref{B}).   For $a \in \F_q$ and $b\in\F_q^*$, we have 
	\begin{align*}
		\left| \Nm_M(a, b) - 
		\left( \frac{|M^*|}{q(q-1)}
		+\frac{(-1)^{\sum_i d_i}q^{\frac{n-m}{2}}}{q} \right) \right|
		\leq (\sum_i d_in_i -1)q^{\frac{n-2}{2}}.
	\end{align*}

	\end{theorem} 
This theorem is proved in \cite{LW} in the case when $M=B$ is commutative, see also \cite{Ka, MW10} for earlier results in the 
non-split case $M=B= \F_{q^n}$ and in the geometric split case $M=B = \F_q^n$. 	

\section{Proof of Theorem \ref{thm}}

We now apply Lemma \ref{lem} to prove the estimate for the higher Kloosterman sum in Theorem \ref{thm}. 
Let $\psi: \F_q \longrightarrow \C^*$ be a fixed non-trivial additive character. Let $f(t) \in \F_q[t]$ be a polynomial of degree $r\geq 1$ with $(r, p)=1$. 
For $b\in \F_q^*$,  recall that the higher Kloosterman sum over $M$ associated to $f(t)$ is defined by
$$K_{M,f}(b, \psi): = \sum_{x \in M^*, \ {\rm N}_M(x)=b} \psi\left(f({\rm Tr}_M(x))\right).$$
By definition, we can rewrite 
$$K_{M, f}(b, \psi)= \sum_{a\in \F_q} \Nm_M(a,b)\psi(f(a)).$$
Recall the notation $m=\sum_i d_in_i$. 
By Lemma \ref{lem}, one computes
\begin{align}\label{kl0}
K_{M, f}(b, \psi) &= \sum_{a\in \F_q} \left( \frac{|M^*|}{q(q-1)}+ \frac{ (-1)^{m - \sum_i d_i}q^{\frac{n-m}{2}}}{q(q-1)} S_m(a, b)\right) \psi(f(a)) \nonumber\\
	&= \frac{|M^*|}{q(q-1)} \sum_{a\in \F_q}\psi(f(a)) + \frac{(-1)^{m - \sum_i d_i}q^{\frac{n-m}{2}}}{q(q-1)} \sum_{v \in \F_q^*}\sum_{a \in \F_q} \psi( f(a)-av) \sum_{\chi} \bar{\chi}(bv^m)
	G(\chi)^{m} \nonumber \\
	&= \frac{|M^*|}{q(q-1)} \sum_{a\in \F_q}\psi(f(a)) + \frac{(-1)^{m - \sum_i d_i}q^{\frac{n-m}{2}}}{q(q-1)} S_{m,f}(b), 
		\end{align}
		where
\begin{align}\label{kl1}
S_{m,f}(b) &= \sum_{v \in \F_q^*}\sum_{a \in \F_q} \psi( f(a)-av) \sum_{\chi} \bar{\chi}(bv^m)
	G(\chi)^{m} \nonumber \\
	&= \sum_{v \in \F_q^*}\sum_{a \in \F_q} \psi( f(a)-av)\sum_{y_1, \cdots, y_m \in \F_q^*}  \sum_{\chi} \chi(\frac{y_1\cdots y_m}{bv^m}) \psi(y_1+\cdots +y_m) \nonumber \\ 	
&= (q-1)\sum_{v \in \F_q^*}\sum_{a \in \F_q} \psi( f(a)-av)\sum_{y_i \in \F_q^*, \ y_1\cdots y_m = bv^m}  
\psi(y_1+\cdots +y_m). 
\end{align}
This exponential sum can be studied via the toric technique. Here, we use $\ell$-adic methods directly as it is easier to handle the trivial factors. 

First, we assume that $r=1$. In this case, $f(t)=ct + f(0)$ for some $c \in \F_q^*$ and $f(a)-av = (c-v)a+ f(0)$. The above sum over $a$ is zero if $v\not=c$, and thus  
\begin{align}\label{kl2}
 S_{m,f}(b)        &= q(q-1)\psi(f(0)) \sum_{y_1 \cdots y_m =bc^m} \psi(y_1+\cdots +y_m)  \\ 
 &=q(q-1)\psi(f(0))K_{\F_q^m}(bc^m, \psi),  
\end{align}
where $K_{\F_q^m}(bc^m, \psi)$ denotes the classical hyper-Kloosterman sum in (\ref{klsum}). 
Since $\sum_{a\in \F_q}\psi(f(a)) = 0$ when $\deg(f)=r=1$, it follows that 
\begin{align}\label{r=1}
K_{M, f}(b, \psi) &= (-1)^{m - \sum_i d_i} \psi(f(0))q^{\frac{n-m}{2}}K_{\F_q^m}(bc^m, \psi).
\end{align}
This proves that 
$$
 |K_{M, f}(b, \psi)| \leq mq^{\frac{n-1}{2}}.
$$
The last inequality agrees with the inequality in Theorem \ref{thm} for $r=1$.  			
In the special case that $M=B$ is \'{e}tale and $f(t)=t$, thus all $d_i=1$ and $n=m$, equation (\ref{r=1}) becomes 
\begin{align}
K_{B, t}(b, \psi) &= (-1)^{n -k}K_{\F_q^n}(b, \psi), 
\end{align}
which is the reduction formula in  Deligne \cite[Equation 7.2.5]{De}. 
In the special case $M= M_{d\times d}(\F_q$) (thus, $k=1$, $n_1=1$, $n=d^2$ and $m=d$) and $f(t)=t$, we obtain 
$$K_{M, t}(b, \psi) = \sum_{x \in {\rm GL}_d(\F_q), \ {\rm det}(x)=b} \psi({\rm Tr}(x)) =q^{{d\choose 2}}K_{\F_q^d}(b, \psi),$$
recovering the reduction formula in Kim \cite[Corollary 5.2]{Ki}.  

Now, we assume that $r\geq 2$ with $(r, p)=1$.  Let $\mathcal{F}= [v \rightarrow bv^m]^* {\rm Kl}_m(\psi)$ be the pullback under the map $v \rightarrow bv^m$ of the standard Kloosterman sheaf ${\rm Kl}_m(\psi)$ on $\mathbb G_m$. Let $\mathcal{G}$ be the Airy sheaf of Katz \cite{K} 
on $\mathbb A^1$ associated to the family of one variable polynomial $f(t) -vt \in \F_q[t]$ parameterized by $v$.  
For $v \in \F_q^*$, the trace formula for 
Kloosterman sheaf and Airy sheaf gives 
\begin{align*}
\sum_{a \in \F_q} \psi( f(a)-av)&=-{\rm Tr}(F_{v, \F_q}| \mathcal{G}),   \nonumber \\
 \sum_{y_i \in \F_q^*, \ y_1\cdots y_m = bv^m} 
\psi(y_1+\cdots +y_m) &= (-1)^{m-1}{\rm Tr}(F_{v, \F_q}| \mathcal{F}). 
\end{align*}
By (\ref{kl1}), and the $\ell$-adic trace formula, we obtain 
\begin{align}\label{kl3}
S_{m, f}(b)&=(-1)^{m}(q-1)\sum_{ v \in \F_q^*} {\rm Tr}(F_{v, \F_q}| \mathcal{F} \otimes \mathcal{G}) \nonumber \\ 
&=(-1)^{m}(q-1)\sum_{i=0}^2 (-1)^i {\rm Tr}(F_{\F_q}| H_c^i( \mathbb{G}_m \otimes \bar{\F}_q, \mathcal{F}\otimes \mathcal{G})).
\end{align}

By Deligne \cite[Theorem 7.8]{De}, $\mathcal{F}$ is a lisse sheaf on $\mathbb G_m$ of rank $m\geq 2$, geometrically Lie-irreducible, 
pure of weight $m-1$, tame at $0$, slopes at $\infty$ all bounded by $1$. On the other hand, by Katz \cite[Theorem 17]{K}, the Airy sheaf $\mathcal{G}$ on  $\mathbb A^1$ is lisse of rank $r-1$, pure of weight $1$, slopes at $\infty$ all equal to $r/(r-1) >1$. It follows that $\mathcal{F}\otimes \mathcal{G}$ is a lisse sheaf on $\mathbb G_m$ of rank $m(r-1)$, no geometrically trivial quotient, pure of weight $m$, tame at $0$, all slopes at $\infty$ equal to $\max(1, r/(r-1)) = r/(r-1)$. This implies that $H_c^0 =H_c^2 = 0$. We obtain 

\begin{align}\label{kl4}
S_{m, f}(b)
&=(-1)^{m+1}(q-1) {\rm Tr}(F_{\F_q}| H_c^1( \mathbb{G}_m \otimes \bar{\F}_q, \mathcal{F}\otimes \mathcal{G})).
\end{align}
By Deligne's theorem, 
 $H_c^1$ is mixed of weight $\leq m+1$.   Furthermore,  as $H_c^0 = H_c^2=0$, 
 $${\rm dim}(H_c^1) = -\chi_c = {\rm Swan}_{\infty}(\mathcal{F}\otimes \mathcal{G}) = m(r-1)\frac{r}{r-1}=mr.$$
 It follows  that 
  $$|S_{m, f}(b) |  \leq mr (q-1)q^{\frac{m+1}{2}}.$$ 
 This is slightly weaker than what we need.  
  
  The cohomology $H_c^1( \mathbb{G}_m \otimes \bar{\F}_q, \mathcal{F}\otimes \mathcal{G})$ 
  is in general not pure because of ramification at $0$.  We now compute the trivial factor at $0$ of the Frobenius action on $H_c^1( \mathbb{G}_m \otimes \bar{\F}_q, \mathcal{F}\otimes \mathcal{G})$. By \cite[Theorem 7.8]{De}, the inertia group $I_0$ at $s=0$ 
 acts unipotently on $\mathcal{F}$ with a single Jordan block, and the geometric Frobenius $F_{\F_q}$ acts trivially  on the rank $1$ invariant $\mathcal{F}^{I_0}$.  Since the Airy sheaf $ \mathcal{G}$ is unramified at $0$, the inertia group $I_0$ acts 
 trivially on  $ \mathcal{G}$.  Let $\mathcal{G}_0$ be the stalk of $\mathcal{G}$ at $0$.  
 Since $\mathcal{F}\otimes \mathcal{G}$ has no geometric global sections, the localization sequence injects its inertia invariants at zero into $H_c^1$. One computes that the trivial factor at $s=0$ is 
 $$\det(I -F_{\F_q}T| (\mathcal{F}\otimes \mathcal{G})^{I_0}) = \det(I -F_{\F_q}T| \mathcal{G}_0) = 1 + (\sum_{a\in \F_q}\psi(f(a)))T +\cdots,$$ 
 which is a polynomial in $T$ of degree $r-1$, pure of weight $1$. We can then write 
 $$\det( I -F_{\F_q}T |H_c^1( \mathbb{G}_m \otimes \bar{\F}_q, \mathcal{F}\otimes \mathcal{G}))= \det(I -F_{\F_q}T| \mathcal{G}_0) P(T),$$  
 where $P(T)$ is a polynomial of degree $rm-r+1$, weight at most $m+1$. It follows that 
 $$\left| {\rm Tr}(F_{\F_q}| H_c^1( \mathbb{G}_m \otimes \bar{\F}_q, \mathcal{F}\otimes \mathcal{G})) + \sum_{a\in \F_q}\psi(f(a))\right| 
 \leq (rm-r+1)q^{\frac{m+1}{2}}. $$
 This and equation (\ref{kl4}) imply the estimate 
 $$\left| S_{m, f}(b)  -(-1)^{m}(q-1) \sum_{a\in \F_q}\psi(f(a))\right| 
 \leq (rm-r+1)(q-1)q^{\frac{m+1}{2}}. $$
 This together with (\ref{kl0}) gives the following theorem which is Theorem \ref{thm}. 
  
 \begin{theorem}\label{thm10}
Let $M$ be a finite semi-simple algebra over $\F_q$ of degree $n=\sum_{i=1}^k d_i^2 n_i \geq 2$ as given in equation (\ref{B}) with $m=\sum_{i=1}^k d_in_i$.   
Let $f(t) \in \F_q[t]$ be a polynomial of degree $r\geq 1$ with $(r, p)=1$. 
For $b\in\F_q^*$, we have 
	
\begin{align}\label{kl}
\left|K_{M, f}(b, \psi) - \frac{|M^*| +(-1)^{\sum_i d_i}(q-1)q^{\frac{n-m}{2}}}{q(q-1)} \sum_{a\in \F_q}\psi(f(a)) \right|   \leq\left( rm -r +1\right) q^{\frac{n-1}{2}}. 
		\end{align}
		\end{theorem}

\begin{cor}\label{cor2} 
Let $M$ be a finite semi-simple algebra over $\F_q$ of degree $n=\sum_{i=1}^k d_i^2 n_i \geq 2$ as given in equation (\ref{B}).   
 Let $f(t) \in \F_q[t]$ be a permutation polynomial of degree $r\geq 1$ with $(r, p)=1$.  For $b\in\F_q^*$, we have 
	\begin{align}
		\left| K_{M, f}(b, \psi) \right| 
		\leq \left( (\sum_{i=1}^k d_in_i )r -r +1 \right)q^{\frac{n-1}{2}}.
	\end{align}

\end{cor}

\section{Further questions}

In this section, we discuss several related further problems beginning with the twisted or exotic version of higher 
Kloosterman sum over the semi-simple algebra $M$. 

\begin{definition}\label{def-ex} 
Given a multiplicative character $\chi: M^* \longrightarrow \C^*$,  the twisted or exotic 
version of higher Kloosterman sum over $M$ is defined by: 
$$K_{M,f}(b, \chi, \psi): = \sum_{x \in M^*, \ {\rm N}_M(x)=b} \chi(x) \psi\left(f({\rm Tr}_M(x))\right),$$
where $f(t) \in \nfq[t]$ and $b\in \nfq^*$ as before. In the case $\chi=1$, this reduces to the sum $K_{M,f}(b, \psi)$ we studied in the previous section. 
\end{definition}

\begin{problem} Extend Theorem \ref{thm10} to this exotic situation with arbitrary $\chi$.  
\end{problem} 

If $\chi$ factors through the norm map ${\rm N}_M$, namely, $\chi =\chi_0 \circ {\rm N}_M$, where $\chi_0: \nfq^* \longrightarrow \C^*$ 
is a multiplicative character of $\nfq^*$, then by definition, one checks  
$$K_{M,f}(b, \chi, \psi) = \chi_0(b) K_{M,f}(b, \psi),$$
reducing to the untwisted case in Theorem \ref{thm10}. So, the new situation is when $\chi$ does not factor through the norm map ${\rm N}_M$. 

\begin{remark}
If $f(t)=t$, the exotic sum $K_{M,t}(b, \chi, \psi)$ becomes
$$K_{M}(b, \chi, \psi): = \sum_{x \in M^*, \ {\rm N}_M(x)=b} \chi(x) \psi\left({\rm Tr}_M(x)\right).$$
When $M=B$ is commutative, this sum $K_{B}(b, \chi, \psi)$ is the exotic Kloosterman sum introduced by Deligne \cite[Remark 7.18]{De} 
and studied by Katz \cite[Theorem 8.8.5]{K0}. When $M$ is non-commutative, the exotic Kloosterman sum 
$K_{M}(b, \chi, \psi)$ is studied in Wan-Zhang \cite{WZ}.  In Definition \ref{def-ex}, we only considered abelian character 
$\chi$ of $M^*$. One could also consider non-abelian character $\chi_{\rho} ={\rm Tr}\circ \rho$ attached to a linear representation $\rho:  M^* \longrightarrow {\rm GL}_d(\mathbb{C})$, 
which would seem to be more difficult.  
\end{remark}

Now, we discuss the connection to hypergeometric sums 
over reductive groups. We will restrict to the untwisted case. Taking $M = M_{d}(\F_q)$ with $d>1$ and $b=1$ in Theorem \ref{thm10}, we obtain 

\begin{cor}\label{cor3}
Let $f(t) = \sum_{i=1}^r a_i t^i \in \F_q[t]$ be a polynomial of degree $r\geq 1$ with $(r, p)=1$.  
For each $1\leq i \leq r$, let $\rho_i$ denote the 
$i$-th tensor power of the standard representation $\rho$ of ${\rm SL}_d(\F_q)$ with $d>1$. Then, we have  
\begin{align}
		\left| \sum_{x \in {\rm SL}_d(\F_q)} \psi(\sum_{i=1}^r a_i {\rm Tr}(\rho_i(x))) 
		- \frac{|{\rm SL}_d(\F_q)| +(-1)^dq^{\frac{d^2-d}{2}} }{q} \sum_{a\in \F_q}\psi(f(a))		\right| 
		\leq (dr -r +1) q^{\frac{d^2-1}{2}}.
	\end{align}
If in addition, $f(t) = \sum_{i=1}^r a_i t^i \in \F_q[t]$ is a permutation polynomial over $\F_q$, then we have 		  
	\begin{align}
		\left| \sum_{x \in {\rm SL}_d(\F_q)} \psi(\sum_{i=1}^r a_i{\rm Tr}(\rho_i(x))) \right| 
		\leq (dr -r +1)q^{\frac{d^2-1}{2}}.
	\end{align}
	\end{cor}

\begin{remark}
More generally, one can consider the following $r$ term hypergeometric exponential sum over the reductive group ${\rm SL}_d$ associated to the 
representations $\rho_i$ ($1\leq i \leq r$): 
\begin{align}
S(A_1, \cdots, A_r): = \sum_{x \in {\rm SL}_d(\F_q)} \psi(\sum_{i=1}^r {\rm Tr}(A_i\rho_i(x))), 
\end{align}
where $A_i$ is an endomorphism of the representation space for $\rho_i$. 
The general result in Fu-Li \cite{FL} implies that for a generic (more precisely non-degenerate) choice of the endomorphisms $A_i$, the sum $S(A_1, \cdots, A_r)$ sometimes has  
the expected square root cancellation. However, no explicit example is given in \cite{FL}. 
In the case $A_i = a_i$ (the scalar multiplication by $a_i$), Corollary \ref{cor3} gives explicit examples of hypergeometric sums over ${\rm SL}_d(\F_q)$ 
with explicit square root cancellation. Because of the possible existence of a non-trivial main term,  Corollary \ref{cor3} cannot be deduced from the 
general results in \cite{FL}, as there will be non-trivial cohomology beyond the middle dimension for $r>1$.  Note that as $f(t)$ varies over the affine space $\mathbb{A}^{r}=\{(a_1, \cdots, a_r)\}$ over $\F_q$, we do not get the square root cancellation 
for the above exponential sum for any $f$ with degree $r>1$ and all $\F_q$, but only for some special $f$ and $\F_q$ such as those $f$ which induces a permutation of $\F_q$. A two term hypergeometric sum over ${\rm GL}_d$, called matrix Kloosterman sum, is studied in 
Erd\'elyi-T\'oth \cite{ET} and Erd\'elyi-Sawin-T\'oth \cite{EST}. More generally, an $r$ term matrix Kloosterman sum is studied in Zelingher \cite{Ze}. 
Its extension to semi-simple algebra $M$ is the following product-trace Kloosterman sum defined by 
$$K(M, r, x, \psi):=\sum_{g_1\cdots g_r=x, ~g_i \in M^*} \psi ({\rm Tr}_M(g_1+\cdots +g_r)),$$
where $r\geq 2$ and $x \in M^*$. This sum has the expected square root cancellation if $x$ is a regular element in $M^*$, 
see \cite[Corollary 4.9]{Ze} and \cite[Theorem 4.2]{Wa2}.

\end{remark}

\begin{remark}
It would be interesting to prove similar explicit results when the $i$-th tensor power $\rho_i$ is replaced by the $i$-th symmetric power 
(or $i$-th Adams operation) of the standard representation $\rho$. In the Adams operation case, the general problem over the semi-simple algebra $M$ is to estimate the sum 

\begin{align}
A_{M,f}(b, \psi): = \sum_{x \in M^*, \ {\rm N}_M(x)=b} \psi\left({\rm Tr}_M(f(x))\right),
\end{align}
where $f(t)\in \F_q[t]$ is a polynomial of degree $r$ not divisible by $p$. Note that now $f$ is inside the trace map ${\rm Tr}_M$, not outside of ${\rm Tr}_M$ as in the definition of $K_{M,f}(b, \psi)$. 
In the commutative split case $M=B =\F_q^n$,  applying the work of 
Adolphson-Sperber \cite[Theorem 4.2]{AS} and Denef-Loeser \cite{DL} for toric exponential sums, one obtains the following estimate:

\begin{align}\left| A_{\F_q^n,f}(b, \psi): = \sum_{x_1, \cdots, x_n \in {\F_q^*}, \ x_1\cdots x_n=b} \psi\left(f(x_1)+\cdots + f(x_n)\right)\right| \leq nr^{n-1} q^{\frac{n-1}{2}}, \ b\in \F_q^*.
\end{align}
In the general commutative finite \'etale case $M=B= \prod_{i=1}^k \F_{q^{n_i}}$ of degree $n=\sum n_i$, the same estimate 
\begin{align}\label{es1}
\left| A_{B,f}(b, \psi): = \sum_{x \in B^*, \ {\rm Norm}(x)=b} \psi\left({\rm Tr}_B(f(x))\right)\right| \leq n r^{n-1} q^{\frac{n-1}{2}}, \ b \in \F_q^*
\end{align}
holds by a twisting or Weil descend argument. This is already a good estimate. 
For non-commutative 
semi-simple algebra $B$ over $\F_q$, we do not know a good estimate for $A_{B,f}(b, \psi)$.

\end{remark}

\begin{remark}
Another direction is to study the following higher inverted Kloosterman sum over semi-simple algebra $M$: 
\begin{align}
IK_{M, f}(b, \psi): =  \sum_{x \in M^*, \ f({\rm Tr}_M(x))\not=0, \ {\rm N}_M(x)=b} \psi(\frac{1}{f({\rm Tr}_M(x))}).
\end{align}
In the case when $M=B$ is commutative and $f(t)=t$,  this inverted Kloosterman sum and its exotic version are introduced in Katz \cite{Kb},  and studied in detail in 
\cite{LW0} and \cite{FW}.  The question is to extend these results from commutative case with $f(t)=t$ to non-commutative case with 
general higher degree polynomial $f(t)$. More generally, one can also try to replace the polynomial $f(t)$ and the reciprocal polynomial $1/f(t)$ by a rational 
function $g(t)$.  The method of the present paper still works in principle, but one first needs to extend Katz \cite[Theorem 17]{K} on Airy sheaf from polynomials to rational functions. Needless to say, there will be non-trivial main terms in such sums, and thus there will be non-trivial cohomology beyond the 
middle dimension. 
\end{remark}

\section{Non-vanishing exponential sums}

In this section, we discuss when the one variable exponential sum 
$$S_q(f): = \sum_{t \in \nfq} \psi (f(t))$$ 
is non-zero, where $\psi$ is a fixed non-trivial additive character of $\nfq$ and $f(t)\in \nfq[t]$ is a polynomial 
of positive degree $r$. For simplicity of later exposition, in this section, without loss of generality, we will take 
$$ \psi = \psi_0 \circ {\rm Tr}_{\F_q /\F_p} = \psi_0 \circ {\rm Tr}$$ 
to be a standard additive character of $\F_q$, where $\psi_0$ is a fixed non-trivial additive character of $\F_p$ and 
${\rm Tr}$ denotes ${\rm Tr}_{\F_q /\F_p}$ in this section.  

One simple relation between permutation polynomials and vanishing exponential sums over $\nfq$ is given below. 

\begin{proposition} A polynomial $f(t) \in \nfq [t]$ is a permutation polynomial over $\nfq$ if and only if $S_q(af)=0$ for all $a\in \nfq^*$. 
\end{proposition}

\begin{proof} If $f(t)$ is a permutation polynomial over $\nfq$, it is clear that $S_q(af)=S_q(t)=0$ for all $a\in \nfq^*$. Conversely, assume that 
$S_q(af)=0$ for all $a\in \nfq^*$. Then, for all $b\in \nfq$, the number of roots of $f(t)-b=0$ in $\nfq$ is given by 
$$\frac{1}{q} \sum_{a,  t \in \nfq} \psi (a(f(t)-b)) = 1 + \frac{1}{q}\sum_{a\in \nfq^*} \psi(-ab) S_q(af) =1.$$
That is, $f(t)$ is a permutation polynomial over $\nfq$. 
\end{proof}

A natural problem is to classify all polynomials $f(t) \in \nfq[t]$ such that $S_q(f)=0$, or just to count the number of such polynomials. 
Write $\zeta_p =\psi_0(1)$, which is a primitive $p$-th root of unity. 
Then, we have the expansion 
$$S_q(f) = \sum_{k=0}^{p-1}\# \{ t\in\nfq \ |\ {\rm Tr}(f(t))=k \} \zeta_p^k.$$
The minimal polynomial of $\zeta_p$ over $\mathbb{Q}$ is $x^{p-1}+ \cdots + x+1$. It follows that $S_q(f)=0$ if and only if 
$\# \{ t\in\nfq \ |\ {\rm Tr}(f(t))=k \} = \frac{q}{p}$ for all $k\in \nfp$. Let 
$$L_k = {\rm Tr}^{-1}(k) = \{a \in \nfq ~|~ {\rm Tr}(a)=k \},$$which is a subset of $\nfq$ with 
cardinality $q/p$. By Lagrange interpolation, any polynomial $f\in \nfq[t]$ of degree at most $q-1$ satisfying $S_q(f)=0$ 
can be uniquely identified with a map of the following type: 
$$f: \ \nfq \longrightarrow \nfq = \bigcup_{k \in \nfp} L_k, \ \ |f^{-1}(L_k)| = \frac{q}{p}, \ \ \ k \in \nfp.$$
A combinatorial argument then gives
$$\# \{ f \in \nfq [t]| \deg(f) \leq q-1, S_q(f)=0\} = {q \choose \frac{q}{p}, \cdots, \frac{q}{p}} (\frac{q}{p})^q.$$
This is significantly more than the number $q!$ of permutation polynomials over $\nfq$ if $q$ is a proper power of $p$. 
If $q=p$, they are the same. For $0\leq r \leq q-1$, let $M_r(q)$ denote the number of $f \in \nfq[t]$ such that $\deg(f)=r$ and 
$S_q(f)=0$. As noted above, we have 
\begin{equation}
m(1, q): = \sum_{r=0}^{q-1} M_r(q) = {q \choose \frac{q}{p}, \cdots, \frac{q}{p}} (\frac{q}{p})^q.
\end{equation}
It would be interesting to determine the second moment 
$$m(2, q): = \sum_{r=0}^{q-1} M_r^2(q).$$
More generally, we raise
\begin{problem} Determine $M_r(q)$ for each individual $r$ or its distribution as $r$ varies. 
\end{problem}

This is expected to be a difficult problem in general. Clearly, $M_0(q)=0$ and $M_1(q)=q(q-1)$. 
If $1< r$ and $r$ divides $p-1$, then a $p$-adic valuation argument shows that $M_r(q)=0$. If $(r, q-1)>1$, then it is conjectured that 
$M_r(q)=0$ if $p$ is sufficiently large compared to $r$. 

For $\alpha \in \nfq$, let $k ={\rm Tr}(\alpha)$. Then, $\{a \in \nfq ~|~ {\rm Tr}(a)=k \} = \{ \alpha + b, b \in \nfq ~|~ {\rm Tr}(b)=0\}$. 
It follows that 
$$\# \{ t\in\nfq \ |\ {\rm Tr}(f(t))=k \}  = \# \{ t \in \nfq ~ | ~{\rm Tr}(f(t) -\alpha)=0\}.$$
This gives 
\begin{proposition}
An equivalent condition for $S_q(f)=0$ is that 
the affine Artin-Schreier curve $y^p -y +\alpha = f(t)$ has exactly $q$ rational points over $\nfq$ for all $\alpha \in \nfq$. 
\end{proposition}

Alternatively, it is clear that $S_q(f)=0$ if and only if $S_q(cf)=0$ for all $c \in \nfp^*$ (the multiplicative group of the prime subfield of $\F_q$). 
Since 
\begin{align*}
\sum_{c \in \nfp^*} | S_q(cf(t))|^2 &= \sum_{c \in \nfp} \sum_{x, y \in \nfq} \psi (c(f(x)-f(y))) - q^2 \\
&= \#\{ (x, y, z) \in \F_q^3: z^p-z=f(x)-f(y)\}-q^2,
\end{align*}
where in the last step, we used our convention that $\psi = \psi_0 \circ {\rm Tr}_{\F_q /\F_p}$ is a standard additive character of $\F_q$.  
This gives 

\begin{proposition}
Let $f(t) \in \nfq[t]$. Then, $S_q(f)=0$ if and only if  the Artin-Schreier surface $z^p-z=f(x)-f(y)$ 
has exactly $q^2$ rational points over $\nfq$.
\end{proposition}

\begin{example}
 For monomial $f(t)= at^r$ ($a\in \nfq^*$), it is easy to see that $S_q(at^r)=0$ 
if and only if $(r, q-1)=1$.  
In fact, if $(r, q-1)>1$, then as algebraic integers, 
$$S_q(at^r) = 1 + (r, q-1) \sum_{x \in (\nfq^*)^{(r, q-1)}} \psi(ax) \equiv 1 \mod (r, q-1).$$
It follows that $S_q(at^r) \not=0$. 
\end{example}

\begin{example}
 The Dickson polynomial 
 $$D_r(t, a): = \left( \frac{t +\sqrt{t^2-4a}}{2}\right)^r + \left( \frac{t -\sqrt{t^2-4a}}{2}\right)^r \in \F_q[t] , \ \ a \in \nfq^*$$ satisfies 
 \[ D_0 =2, D_1 =t,  D_r= tD_{r-1} - aD_{r-2}.
 \]
 It is clear that $S_q(D_r(t, a))=0$ if $(r, q^2-1)=1$, as in this case, $D_r(t, a)$ is a permutation polynomial over $\F_q$. 
The converse is not true. That is, the condition $(r, q^2-1)>1$ does not imply that $S_q(D_r(t, a))\not= 0$. In fact, even the 
stronger condition $r|(q-1)$ does not imply that $S_q(D_r(t, a))\not= 0$. The following is an explicit counter-example with $q=25$ and $r=3$.  
Take $\mathbb{F}_{25} = \mathbb{F}_5(\theta)$ with $\theta^2=2$. Then, 
\[ D_3(t, 2\theta) = t^3 -\theta t, \  3 |(25-1). 
\]
Write $t = u+ v\theta$ with $u, v \in \mathbb{F}_5$, half the absolute trace of $t^3-\theta t$ is $u^3 +uv^2 -2v$. For $u=0$, the exponential sum 
over $v$ vanishes. For $u\not=0$, use $u^3 = u^{-1}$ in $\mathbb{F}_5$ to rewrite the expression as $u(v-u^{-1})^2$. The remaining exponential sum 
is 
\[  \sum_{u \in \mathbb{F}_5^*}\sum_{v\in \mathbb{F}_5} \psi_0(2uv^2) = 4  + 4(-1) =0.
\]
\end{example}

\begin{remark}
It would be interesting to classify 
all polynomials $f(t) \in \nfq[t]$ such that $S_{q^e}(f)=0$ for infinitely many $e$. In this case, we say that $f(t)$ gives exceptional exponential sums. 
Since $S_q(f(t) +g(t)^p) = S_q(f(t)+g(t))$, we can always reduce all terms in $f$ whose exponents are divisible by $p$. 
Thus, without loss of generality, we can assume that $f(t)$ is Artin-Schreier reduced, in particular, the degree $r$ of $f(t)$ is not divisible by $p$. 
\end{remark}

If $(r, q-1)=1$, then $S_{q^e}(t^r)=0$ for infinitely many $e$. 
One can hope that the converse might be true for the degree $r$ of $f(t)$ if the characteristic $p$ is large  compared to $r$. 
The following is an exponential sum analogue of Lenstra's theorem on the Carlitz-Wan conjecture for permutation polynomials. 

\begin{conjecture} \label{conj}
Let $(r, q-1)>1$. 
If the characteristic $p$ is sufficiently large compared to $r$,  then $S_q(f) \not =0$ for all $f(t) \in \nfq [t]$ of degree $r$. 
\end{conjecture}

Here we give a geometric heuristic which suggests that the conjecture should be true in many cases. Since $p$ is large compared to $r$, we can assume that 
$r$ is not divisible by $p$. 
For $a \in \nfq$, let $N_{a}$ denote the 
number of $\nfq$-rational points on the plane curve $f(x)-f(y)=a$. Let $n_{a}$ denote the number of 
absolutely irreducible components of the curve $f(x)-f(y)=a$ over $\nfq$. By the Weil bound, we have 
$$N_{a} = n_{a} q + E_{a}, \ \ E_{a} = O_r(\sqrt{q}).$$
The plane curve $f(x)-f(y)=a$ over $\nfq$ has at least one non-singular $\nfq$-rational point at infinity, namely, $(1,1,0)$. Bezout's theorem implies that 
the curve $f(x)-f(y)=a$ has at least one absolutely irreducible component over $\nfq$ for every $a \in \nfq$. That is, 
$n_{a} \geq 1$ for all $a \in \nfq$. Now $(r, q-1)>1$ and $r$ is not divisible by $p$, the plane curve $f(x)-f(y)=0$ has at least two 
non-singular $\nfq$-rational points at infinity, namely, $(1, 1, 0)$ and $(1, \zeta, 0)$, where $\zeta$ is a primitive $(r, q-1)$-th root of unity in $\nfq^*$. Removing the simple trivial factor $x-y$, one sees that the plane curve $\frac{f(x)-f(y)}{x-y} =0$ 
still has at least one non-singular $\nfq$-rational point at infinity, namely, $(1, \zeta, 0)$. This and Bezout theorem show that $n_0\geq 2$. 

Assume now that $S_q(f)=0$. Then  the surface $z^p-z = f(x)-f(y)$ has exactly $q^2$ rational points over $\nfq$. It follows that  
$$q^2 = p\sum_{ {\rm Tr}(a)=0} N_{a} = p\sum_{ {\rm Tr}(a)=0} n_{a} q + p
\sum_{ {\rm Tr}(a)=0} E_{a}\geq (q+p)q  + p\sum_{ {\rm Tr}(a)=0} E_{a}.$$
This implies that 
 \begin{equation}\label{5}
 |\sum_{ {\rm Tr}(a)=0} E_{a} |  \geq q.
 \end{equation}
 As $a$ varies in the trace-zero subset $\{ a \in \F_q:  {\rm Tr}(a)=0\}$ with cardinality $q/p$, in many cases, one expects that the error term $E_{a}$ is randomly distributed with at least the exact correlation $E_a = E_{-a}$ obtained by interchanging $x$ and $y$. In such cases, one would expect that  
  $$|\sum_{ {\rm Tr}(a)=0} E_{a} |  = O_r(\sqrt{\frac{q}{p}} \sqrt{q}) = O_r(\frac{q}{\sqrt{p}}).$$
 This is certainly smaller than $q$ if $p$ is sufficiently large compared to $r$, which is a contradiction to inequality (\ref{5}).

\section{Another relation}
In this final section, we now give another relation between vanishing exponential sums and permutation polynomials.

\begin{proposition} \label{prop5}
Let $f(t) \in \nfq[t]$ be a polynomial of degree at most $q-1$. Then, $S_q(f)=0$ if and only if there is a polynomial 
$h(t)\in \nfq[t]$ of degree at most $q-1$ such that $f(t)  +  h(t)^p - h(t)$ is a permutation polynomial over $\nfq$. 
\end{proposition}

\begin{proof}The ``if" part is obvious. We now prove the ``only if" part. 
Assume $S_{q}(f)=0$. 
Recall that for $k \in \nfp$, we defined $L_k = {\rm Tr}^{-1}(k)$, which is a subset of $\nfq$ with 
cardinality $q/p$. Any polynomial $f\in \nfq[t]$ of degree at most $q-1$ satisfying $S_q(f)=0$ 
can be uniquely identified with a map of the following type: 
$$f: \ \nfq \longrightarrow \nfq = \bigcup_{k \in \nfp} L_k, \ \ |f^{-1}(L_k)| = q/p, \ \ \ k \in \nfp.$$
This map $f$ is not a permutation of $\nfq$ in general, as $f: ~f^{-1}(L_k) \longrightarrow L_k$ might not be injective for  some $k \in \nfp$. 
Namely, we only have the inclusion $f(f^{-1}(L_k)) \subseteq L_k$. Now, we can simply modify (non-uniquely) the values of $f$ on each 
subset $f^{-1}(L_k)$ to obtain a map $g: f^{-1}(L_k) \longrightarrow L_k$ which is bijective for all $k \in \nfp$. By Lagrange interpolation, such a map $g$ is then given by a permutation 
polynomial of degree at most $q-1$ and satisfies  ${\rm Tr}(f(t)) = {\rm Tr}(g(t))=k$ for all $t \in f^{-1}(L_k)$ and all $k \in \nfp$. 
In particular, ${\rm Tr}(f(t)) = {\rm Tr}(g(t))$ for all $t \in \mathbb{F}_q$. For each $t\in \nfq$, choose $\alpha_t \in \nfq$ such that $g(t) = f(t) + \alpha_t^p -\alpha_t$. 
Lagrange interpolation again gives a polynomial $h(t) \in \nfq[t]$ of degree at most $q-1$ such that $h(t)= \alpha_t$ for all $t \in \nfq$.  This implies that 
$g(t) = f(t)+ h(t)^p -h(t)$ for all $t \in \mathbb{F}_q$, 
that is, 
$$g(t) = f(t)+ h(t)^p -h(t) \mod (t^q-t),$$
where $g(t) \in \nfq[t]$ is a permutation polynomial of degree at most $q-1$. 
\end{proof}

\begin{remark}
The above proposition suggests a method to construct non-permutation polynomial $f(t)\in \nfq [t]$ 
such that $S_{q}(f)=0$: Start with a permutation polynomial $g(t)\in \mathbb{F}_q[t]$, usually not Artin-Schreier reduced, and add or subtract a 
polynomial of the form $h(t)^p - h(t)$ with $h(t) \in \mathbb{F}_q[t]$. 
The new polynomial $f(t) = g(t) + h(t)^p -h(t)$ or $f(t) = g(t) - (h(t)^p -h(t))$ is sometimes a non-permutation polynomial such that 
$S_q(f)=S_q(g)=0$. 
\end{remark}

\begin{example}   A trivial example is to take $r=pr_1$ where $(r_1, p-1)>1$ and $f(t)= t^{pr_1} - t^{r_1}+t$. Clearly, $S_{p^e}(f(t)) = S_{p^e}(t)=0$ 
for all $e\geq 1$, but $(r, p^e-1)>1$. For $p^e > (pr_1)^4$, Lenstra's theorem implies that the polynomial $f(t)$ is not a permutation polynomial over $\F_{p^e}$. In this example, $p$ is not 
sufficiently large compared to $r=pr_1$, although $q=p^e$ can be as large as one wishes compared to $r$. 
This example shows that Conjecture \ref{conj} is false if one replaces the condition that $p$ is large 
compared to $r$ with the weaker condition that $q$ is large 
compared to $r$.
\end{example}

\begin{example}
A non-trivial example with $(r, p)=1$ is given by Zhiguo Ding. Let $p$ be odd, $q=p^e$ and $a\in \mathbb{F}_q^*$. Let  
$$f_a(t)= t(t^\frac{p-1}{2} -a)^2\in \mathbb{F}_q[t].$$
 This is a polynomial of degree $p$. 
 A classical result of Dickson says that $f_a(t)$ 
is a permutation polynomial over $\mathbb{F}_q$ if and only if $t^{(p-1)/2}-a$ has  
no roots in $\mathbb{F}_q^*$. 
Take any such $a \in \mathbb{F}_q^*$. Then, $S_q(f_a(t))=0$ since $f_a(t)$ is a permutation polynomial. 
Now, 
$$f_a(t)= t^p -2a t^\frac{p+1}{2} + a^2t = t^p-t -2a t^\frac{p+1}{2} + (a^2+1)t=t^p-t+ g_a(t),$$
where $g_a(t) = -2a t^\frac{p+1}{2} + (a^2+1)t$. 
It follows that  $S_q(g_a(t)) = S_q(f_a(t))=0$. Note that $g_a(t)$ is a polynomial of degree 
$(p+1)/2$ which divides $(q-1)=(p^e-1)$ if $e$ is even. This implies that for $e$ even and 
$a\in \mathbb{F}_q^*$ such that $t^{(p-1)/2}-a$ has  
no roots in $\mathbb{F}_q^*$, the polynomial $g_a(t)$ is a polynomial of degree $(p+1)/2$ 
which divides $(q-1)$ and thus $g_a(t)$ is not a permutation polynomial over $\mathbb{F}_q$ by Hermite's criterion.  
\end{example} 
 
 By Proposition \ref{prop5}, Conjecture \ref{conj} is equivalent to the following 
   
 \begin{conjecture}\label{conj2}
Let $(r, q-1)>1$. 
If the characteristic $p$ is sufficiently large compared to $r$,  then for every polynomial $f(t) \in \nfq[t]$ of degree $r$ and every polynomial $h(t) \in \nfq[t]$, the new polynomial 
$f(t) + h(t)^p -h(t)$ is not a permutation polynomial over $\nfq$.  
\end{conjecture}  
    
As mentioned before, the conjecture is clearly true if $q=p$, since in this case, $h(t)^p -h(t)$ is the zero function on $\mathbb{F}_p$. One can then apply 
the results on the easier tame case of the Carlitz-Wan conjecture. 
If $1< r \ |\ (p-1)$, then viewed as a $p$-adic number, the $q$-adic valuation is given by 
$$v_q\left(\sum_{a\in \F_q}\psi(f(a))\right) = \frac{1}{r}, $$
and thus there are 
no polynomials $f(t)$ of degree $r$ with $\sum_{a\in \F_q}\psi(f(a))=0$. In particular, the conjecture is true if $r \ |\ (p-1)$, in which case, one does not need to assume that $p$ is sufficiently large compared to $r$.  If $(r, q-1)=1$, then $\sum_{a\in \F_q}\psi(a^r)=0$ and thus the above condition $(r, q-1)>1$ cannot be dropped. 
By a theorem of Cohen \cite{Co}, an exceptional polynomial over $\nfq$ is always a permutation polynomial over $\nfq$. A weaker 
consequence of Conjecture \ref{conj2} is the following 

\begin{conjecture}\label{conj3}
Let $(r, q-1)>1$. 
If the characteristic $p$ is sufficiently large compared to $r$,  then for every polynomial $f(t) \in \nfq[t]$ of degree $r$ and every polynomial $h(t) \in \nfq[t]$, the new polynomial 
$f(t) + h(t)^p -h(t)$ is not an exceptional polynomial over $\nfq$.  
\end{conjecture}


\begin{thebibliography}{12}
	
	\bibitem{AS}A. Adolphson and S. Sperber, Exponential sums and Newton polyhedra: cohomology and estimates, 
	\emph{Ann. Math.}, 130(1989), 367-406. 
	
	\bibitem{Ba} P. Bacik,  Completing the picture for the Skolem problem on order-4 linear recurrence sequences, 
	\emph{TheoretiCS}, Vol. 4 (2025), Article 28, 1-11.
		
	\bibitem{Co}S.D. Cohen, The distribution of polynomials over finite fields, 	
	\emph{Acta Arith.}, Vol. 17 (1970), 255-271.	
	
	
	\bibitem{CF}S.D. Cohen and M.D. Fried, Lenstra's proof of the Carlitz-Wan conjecture on exceptional 
	polynomials: an elementary version, \emph{Finite Fields \& Appl.}, Vol. 1, 3(1995), 372-375.	
	
	\bibitem{DXZ} Z.G.  Ding, W. Xiong and Q.F.  Zhang, Exceptional extensions of local fields and 
	the Carlitz--Wan conjecture, \emph{Sci China Math.}, 68 (2025), 2815-2826.  	
	
	\bibitem{EI} M.  Eichler, Allgemeine Kongruenz-Klasseneinteilungen der Ideale einfacher Algebren 
	\"uber algebraischen Zahlk\"orpern und ihre L-Reihen. \emph{J. Reine Angew. Math.}, 179 (1937), 227 -- 251.

                                 \bibitem{ET}M. Erd\'elyi and \'A. T\'oth, Matrix Kloosterman sums, 
                                \emph{Algebra and Number Theory}, Vol 18 (2024), no. 12, 2247- 2308.     

                                \bibitem{EST}M. Erd\'elyi, W. Sawin and \'A. T\'oth, The purity locus of matrix Kloosterman sums, 
                                \emph{Trans. Amer. Math. Soc.}, Vol 377 (2024), no. 6, 4117- 4132. 
                                
                                                            
                                
	
		
	
	
	
	\bibitem{DL}J. Denef and F. Loeser, Weights of exponential sums, intersection cohomology, and Newton polyhedra, 
	\emph{Invent. Math.}, Vol. 106, 1(1991), 275-294. 
			
	
	\bibitem{De}P. Deligne, “Cohomologie \'{e}tale (SGA 4 $\frac{1}{2}$)", \emph{Lecture Notes in Mathematics},           
	Vol. 569,  Springer-Verlag, Berlin/Heidelberg/New York, 1977.	
	
	\bibitem{FGS} M. Fried, R. Guralnick and Y. Saxl,  Schur covers and Carlitz’s conjecture, 
	\emph{Israel J. Math}. 82 (1993), 157- 225. 	
	
	\bibitem{FL} L. Fu and X. Li, Hypergeometric $\ell$-adic sheaves over reduction groups,  https://arxiv.org/abs/2411.11215. 	
	
	
	
	\bibitem{FW} L. Fu and D. Wan, The exotic inverted Kloosterman sums, \emph{IMRN}, 12 (2025), 1-22. 
	
	\bibitem{GM} R. Guralnick and P. M\"uller, Exceptional polynomials of affine type,  
	\emph{J.  Algebra}, 194 (1997),  429-454.

	\bibitem{He} T. Helleseth,  Some results about the cross-correlation function between two maximal linear sequences, 
	\emph{Discrete Math.}, 16(3):209–232, 1976.	
	
	\bibitem{HZ} Y. Hu and Z. Zhang, A new proof of Carlitz-Wan conjecture on exceptional polynomials,  arXiv:2602.19779. 		
	\bibitem{K} N. Katz, On the monodromy groups attached to certain family of exponential sums, \emph{Duke Math., J.}, Vol 54 (1987), 
	No. 1, 41-- 56. 
		
	\bibitem{K0} N. Katz, Gauss Sums, Kloosterman Sums, and Monodromy Groups, 
	Princeton University Press, 1988. 
		
	\bibitem{Ka} N. Katz, Estimates for Soto-Andrade sums, \emph{J. Reine
	Angew. Math.}, 438(1993), 143--161. 
	
	\bibitem{Kb} N. Katz, A note on exponential sums, \emph{Finite Fields Appl.} 1 (1995), 395-398. 


	
	
	
	
		
	
	
	\bibitem{Kay}N. Kayal, Recognizing permutation functions in polynomial time, \emph{Electronic Colloquium on 
	Computational Compelxity}, Report No. 8, 2005.  	
	
	
	\bibitem{Ki} D. Kim, Gauss sums for general and special linear groups over a finite field, \emph{Arch. Math.}, 69 (1997),  297--304. 	
	
	\bibitem{La} E. Lamprecht, Struktur und Relationen allgemeiner Gaußcher Summen in endlichen Ringen I, II, 
	\emph{J. Reine Angew. Math.},  197 (1957) 1--48.	
	
	
		
	\bibitem{LW0} X. Lin and D. Wan, On inverted Kloosterman sums over finite fields, \emph{Math. Zeit.} 306 (2024) no. 4, 
	Paper No. 60, 22 pp.
	
	
		
	\bibitem{LW} X. Lin and D. Wan, Counting elements with given trace and norm in étale algebras, \emph{International  J. Number Theory}, 
	Vol. 21, 8 (2025), 1955-1965.
	
					
	
	\bibitem{MW10}M. Moisio and D. Wan, On Katz's bound for the number of elements 
	with given trace and norm, \emph{J. Reine Angew. Math.}, 638(2010), 69-74.	
	
		
	
	
	\bibitem{Wa1}D. Wan,  A generalization of the Carlitz conjecture. In: \emph{Finite fields, Coding theory, 
	and Advances in Communications and Computing}, Proceedings of the international conference held 
	at the University of Nevada, Las Vegas, Nevada, August 7–10, 1991. Lecture Notes in Pure and Appl Math, vol. 141. 
	New York: Marcel Dekker, 1993, 431–432.	
	
	
	\bibitem{Wa2} D. Wan, Norm-trace and Kloosterman sums in finite semi-simple algebras, \emph{Frontiers in Combinatorics and 
	Number Theory}, Vol. 1 (2026), 74-88.  
	
	\bibitem{WZ} D. Wan and D. Zhang, Exotic and inverted Kloosterman sums over semisimple algebras, 
	\emph{Finite Fields Appl.}, Vol. 119 (2027), article 102947. 
			
	\bibitem{Ze} E. Zelingher, On matrix Kloosterman sums and Hall-Littlewood polynomials, 
	\emph{Trans. Amer. Math. Soc.}, 378 (2025), 3597–3623.	
	
	
\end{thebibliography}
\end{document}